\documentclass[10pt,twoside]{article}
\usepackage{graphicx}
\usepackage{amsmath}
\usepackage{amssymb}
\usepackage{Latex-document}

\newcommand{\gl}{\mbox{\upshape GL}}
\renewcommand{\H}{\mbox{\upshape H}}
\newcommand{\gln}{\gl_n}
\newcommand{\glm}{\gl_m}

\newcommand{\glnkv}{\gln (k_v)}
\newcommand{\glna}{\gln (\AA )}
\newcommand{\SO}{\mbox{\upshape SO}}
\newcommand{\Sp}{\mbox{\upshape Sp}}
\newcommand{\Ex}{\mbox{\upshape E}}
\newcommand{\Spin}{\mbox{\upshape Spin}}
\newcommand{\CC}{{\mathbb C}}
\renewcommand{\AA}{{\mathbb A}}
\newcommand{\bs}{\backslash}
\newcommand{\diag}{\mathop{diag}}

\newcommand{\Ind}{\mathop{Ind}}
\newcommand{\LH}{{^L\!\H}}
\newcommand{\veps}{\varepsilon}

\title{\bf Converse Theorems, Functoriality, \vskip -2mm
and Applications to Number Theory }
\author{{\bf J. W. Cogdell}\thanks{Department of Mathematics, Oklahoma
State University, Stillwater, OK 74078, USA. E-mail: cogdell@math.okstate.edu} \quad {\bf I. I.
Piatetski-Shapiro}\thanks{Department of Mathematics, Yale University, New Haven, CT 06520, USA, and School of
Mathematics, Tel Aviv University, Tel Aviv 69978, Israel. E-mail: ilya@math.yale.edu}\vspace*{-0.5cm}}
\date{\vspace{-8mm}}

\begin{document}

\maketitle

\thispagestyle{first} \setcounter{page}{119}

\begin{abstract}

\vskip 3mm

There has been a recent coming together of the Converse Theorem
for $\gln$ and the Langlands-Shahidi method of controlling the
analytic properties of automorphic $L$-functions which has allowed
us to establish a number of new cases of functoriality, or the
lifting of automorphic forms. In this article we would like to
present the current state of the Converse Theorem and outline  the
method one uses to apply the Converse Theorem to obtain liftings.
We will then turn to an exposition of the new liftings and some of
their applications.

\vskip 4.5mm


\end{abstract}

\vskip 12mm

\section{Introduction} \label{section 1}\setzero

\vskip-5mm \hspace{5mm}

Converse Theorems traditionally have provided a way to characterize
Dirichlet series associated to modular forms in terms
of their analytic properties.  Most familiar are the Converse Theorems of
Hecke and Weil. Hecke first proved that $L$-functions associated to
modular forms enjoyed ``nice'' analytic properties and then proved
``Conversely'' that these analytic properties in fact
characterized modular $L$-functions. Weil extended this Converse
Theorem to $L$-functions of modular forms with level.

In their modern formulation, Converse Theorems are stated in terms of
automorphic
representations of $\glna$ instead of modular
forms. Jacquet, Piatetski-Shapiro, and Shalika have proved that the
$L$-functions associated to automorphic representations of
$\glna$ have nice analytic properties via integral
representations similar to those of Hecke. The relevant ``nice''
properties are: analytic continuation,
boundedness in vertical strips, and functional equation.
Converse Theorems in this context
invert these integral representations. They give a criterion for an
irreducible admissible representation $\Pi$ of $\glna$ to be
automorphic and cuspidal in terms of the analytic properties
of Rankin-Selberg convolution
$L$-functions $L(s,\Pi\times\pi')$ of $\Pi$ twisted by cuspidal
representations $\pi'$ of $\glm(\AA)$ of smaller rank
groups.

To use Converse Theorems for applications, proving that certain
objects are automorphic, one must be able to show that certain
$L$-functions are ``nice''.  However, essentially the only way to show
that an $L$-function is nice is to have it associated to an
automorphic form. Hence the most natural applications of
Converse Theorems are to functoriality, or the lifting of automorphic
forms,  to $\gln$. More
explicit number theoretic applications then come as consequences of
these liftings.

Recently there have been several applications of Converse Theorems to
establishing functorialities. These have been possible thanks  to the
recent advances in the Langlands-Shahidi method of analysing the
analytic properties of general automorphic $L$-functions, due to
Shahidi and his collaborators \cite{S3}.
By combining our Converse Theorems with their control of the
analytic properties of $L$-functions many new examples of functorial
liftings to $\gln$ have been established. These are described in
Section 4 below. As one number theoretic consequence of these liftings
Kim and Shahidi have been able to establish the best general estimates over a
number field towards the Ramamujan-Selberg conjectures for $\gl_2$,
which in turn have already had other applications.

\section{Converse Theorems for \boldmath$\gln$} \label{section 2} \setzero

\vskip-5mm \hspace{5mm}

Let $k$ be a global field, $\AA$ its adele ring, and $\psi$ a fixed
non-trivial  (continuous) additive character of $\AA$ which is trivial
on $k$. We will take $n\geq 3$ to be an integer.

To state these Converse Theorems, we begin with an irreducible admissible
representation $\Pi$ of $\gl_n(\AA)$.
It has a decomposition
$\Pi=\otimes'\Pi_v$, where $\Pi_v$ is an irreducible admissible
representation of $\gl_n(k_v)$. By the local theory of Jacquet,
Piatetski-Shapiro, and Shalika \cite{JPSS,JS} to each $\Pi_v$ is associated a
local $L$-function $L(s,\Pi_v)$ and a local $\varepsilon$-factor
$\varepsilon(s,\Pi_v,\psi_v)$. Hence formally we can form
$$
L(s,\Pi)=\prod L(s,\Pi_v) \quad\quad\text{ and }\quad\quad
\varepsilon(s,\Pi,\psi)=\prod \varepsilon(s,\Pi_v,\psi_v).
$$
We will always assume the following two things about $\Pi$:
\begin{enumerate}
\item[(1)] $L(s,\Pi)$ converges in some half plane $Re(s)>>0$,
\item[(2)] the central character $\omega_\Pi$ of $\Pi$ is automorphic, that
is, invariant under $k^\times$.
\end{enumerate}
Under these assumptions, $\varepsilon(s,\Pi,\psi)=\varepsilon(s,\Pi)$
is independent of our choice of $\psi$ \cite{CPS1}.

As in Weil's case, our Converse Theorems will involve twists but now
by cuspidal automorphic representations of $\gl_m(\AA)$ for certain
$m$. For convenience, let us set $\mathcal A(m)$ to be the set of automorphic
representations of $\gl_m(\AA)$, $\mathcal A_0(m)$ the set of (irreducible)
cuspidal automorphic representations of $\gl_m(\AA)$, and
$ \mathcal
T(m)=\bigcup_{d=1}^m \mathcal A_0(d)$. If $S$ is a finite set of places,
we will let $\mathcal T^S(m)$ denote the subset of representations
$\pi\in \mathcal T$ with local components $\pi_v$  unramified at all places
$v\in S$ and let $\mathcal T_S(m)$ denote those $\pi$ which are
unramified for all $v\notin S$.

Let $\pi'=\otimes'\pi'_v$ be a cuspidal
representation of $\gl_m(\AA)$ with $m<n$.  Then again we can formally
define
$$
L(s,\Pi\times \pi')=\prod L(s,\Pi_v\times \pi'_v) \quad\quad\text{ and
}\quad\quad
\varepsilon(s,\Pi\times \pi')=\prod
\varepsilon(s,\Pi_v\times \pi'_v,\psi_v)
$$
since the local factors make sense whether $\Pi$ is automorphic
or not. A consequence of (1) and (2) above and the cuspidality of ${\pi'}$ is
that both $L(s,\Pi\times{\pi'})$ and
$L(s,\widetilde\Pi\times\widetilde{\pi'})$
converge absolutely for $Re(s)>>0$, where $\widetilde\Pi$ and
$\widetilde{\pi'}$ are the contragredient representations, and
that $\varepsilon(s,\Pi\times{\pi'})$ is independent of
the choice of $\psi$.

We say that $L(s,\Pi\times{\pi'})$ is {\it nice} if it satisfies the
same analytic properties it would if $\Pi$ were cuspidal, i.e.,
\begin{enumerate}
\item $L(s,\Pi\times{\pi'})$ and
$L(s,\widetilde\Pi\times\widetilde{\pi'})$ have
continuations to {\it entire} functions of $s$,
\item these entire  continuations are {\it bounded in vertical strips} of
finite width,
\item they satisfy the standard {\it functional equation}
\[
L(s,\Pi\times{\pi'})=\veps(s,\Pi\times{\pi'})
L(1-s,\widetilde\Pi\times\widetilde{\pi'}).
\]
\end{enumerate}

The basic converse theorem for $\gl_n$ is the following.


{\bf Theorem 1.} \cite{CPS2} \it Let $\Pi$ be an irreducible admissible
representation of $\gl_n(\AA)$ as above. Let $S$ be
a finite set of finite places. Suppose that
$L(s,\Pi\times{\pi'})$ is nice for all ${\pi'}\in\mathcal T^S(n-2)$.
Then $\Pi$
is quasi-automorphic in the sense that there is an automorphic
representation $\Pi'$ such that $\Pi_v\simeq\Pi'_v$ for all $v\notin S$.
If $S$ is empty, then in fact $\Pi$ is a cuspidal automorphic
representation of $\glna$.\rm


It is this version of the Converse Theorem that has been used in
conjunction with the Langlands-Shahidi method of controlling
analytic properties of $L$-functions in the new examples of
functoriality explained below.


{\bf Theorem 2.} \cite{CPS1} \it Let $\Pi$ be an irreducible admissible
representation of $\gl_n(\AA)$ as above. Let $S$ be
a non-empty finite set of places,  containing $S_\infty$, such that the
class number of the ring $\mathfrak o_S$ of $S$-integers is one. Suppose that
$L(s,\Pi\times{\pi'})$ is nice for all ${\pi'}\in\mathcal T_S(n-1)$. Then $\Pi$
is quasi-automorphic in the sense that there is an automorphic
representation $\Pi'$ such that $\Pi_v\simeq\Pi'_v$ for all $v\in S$
and all $v\notin S$ such that both $\Pi_v$ and $\Pi'_v$ are unramified.\rm


This version of the Converse
Theorem was specifically designed to investigate functoriality in the
cases where one controls the $L$-functions by means of integral
representations where it is expected to be more difficult to control twists.

The proof of Theorem 1 with $S$ empty and $n-2$ replaced by $n-1$
essentially follows the lead of Hecke, Weil, and
Jacquet-Langlands. It is based on the integral representations of
$L$-functions, Fourier expansions, Mellin inversion, and finally a
use of the weak form of Langlands spectral theory. For Theorems 1
and 2 where we have restricted our twists either by ramification
or rank we must impose certain local conditions to compensate for
our limited twists. For Theorem 1 are a finite number of local
conditions and for Theorem 2 an infinite number of local
conditions. We must then work around these by using results on
generation of congruence subgroups and either weak approximation
(Theorem 1) or strong approximation (Theorem 2).

As for our expectations of what form the Converse
Theorem may take in the future, we refer the reader to the last
section of \cite{CPS2}.

\section{Functoriality via the Converse Theorem} \label{section 3}
\setzero\vskip-5mm \hspace{5mm }

In order to apply these theorems, one must be able to control the
analytic properties of the $L$-function. However the only way we
have of controlling global $L$-functions is to associate them to
automorphic forms or representations. A minute's thought will then
convince one that the primary application of these results will be
to the lifting of automorphic representations from some group $\H$
to $\gln$.

Suppose that $\H$ is a reductive group over  $k$. For simplicity of
exposition we
will assume throughout that $\H$ is split and deal only with the
connected component of its $L$-group, which we will (by abuse of
notation) denote by $\LH$ \cite{B}.
Let  $\pi=\otimes'\pi_v$ be a cuspidal automorphic
representation of $\H$ and $\rho$ a complex representation of
$\LH$. To this situation Langlands has associated an $L$-function
$L(s,\pi,\rho)$ \cite{B}. Let us assume that
$\rho$ maps  $\LH$ to $\gln(\CC)$.  Then by Langlands' general
Principle of Functoriality  to $\pi$ should be
associated an automorphic representation $\Pi$ of $\glna$
satisfying $L(s,\Pi)=L(s,\pi,\rho)$,
$\varepsilon(s,\Pi)=\varepsilon(s,\pi,\rho)$, with similar equalities locally
and for the twisted versions \cite{B}.
Using the Converse Theorem to
establish such liftings involves three steps:
 construction of a candidate lift,  verification that the twisted
$L$-functions are ``nice'', and  application of the appropriate
Converse Theorem.


1. {\it Construction of a candidate lift}: We  construct a
candidate lift $\Pi=\otimes'\Pi_v$ on $\glna$ place by place. We
can see what $\Pi_v$ should be at almost all places. Since we have
the arithmetic Langlands (or Hecke-Frobenius) parameterization of
representations $\pi_v$ of $\H(k_v)$ for all archimedean places
and those non-archimedean places where the representations are
unramified \cite{B}, we can use these to associate to $\pi_v$ and
the map $\rho_v: ^L\!\H_v\rightarrow \LH\rightarrow \gln(\CC)$  a
representation $\Pi_v$ of $\glnkv$. This correspondence preserves
local $L$- and $\varepsilon$-factors
\[
L(s,\Pi_v)=L(s,\pi_v,\rho_v) \quad\quad
\text{and}\quad\quad \varepsilon(s,\Pi_v,\psi_v)=\varepsilon(s,\pi_v,\rho_v,\psi_v)
\]
along with the twisted versions.
If $\H$ happens to be $\glm$ or a related group  then
we in principle know how to associate the representation $\Pi_v$ at
all places now that the local Langlands conjecture has been solved for
$\glm$. For other
situations, we may not know what $\Pi_v$ should be at the ramified
places. We will return to this difficulty momentarily and show how one
can work around this with the use of a highly ramified
twist. But for now,
let us assume we can finesse this local problem and arrive at a
global representation $\Pi=\otimes'\Pi_v$ such that
\[
L(s,\Pi)=\prod L(s,\Pi_v)=\prod L(s,\pi_v,\rho_v)=L(s,\pi,\rho)
\]
and similarly $\varepsilon(s,\Pi)=\varepsilon(s,\pi,\rho)$
with similar equalities for the twisted versions.
$\Pi$ should then be the Langlands lifting
of $\pi$ to $\glna$ associated to $\rho$.


2. {\it Analytic properties of global $L$-functions}: For
simplicity of exposition, let us now assume that $\rho$ is simply
a standard embedding of $^L\!\H$ into $\gln(\CC)$, such as will be
the case if we consider $\H$ to be a split classical group, so
that $L(s,\pi,\rho)=L(s,\pi)$ is the standard $L$-function of
$\pi$. We have our candidate $\Pi$ for the lift of $\pi$ to $\gln$
from above. To be able to assert that the $\Pi$ which we
constructed place by place is automorphic, we will apply a
Converse Theorem. To do so we must control the twisted
$L$-functions $L(s,\Pi\times\pi')=L(s,\pi\times\pi')$ for
$\pi'\in\mathcal T$ with an appropriate twisting set $\mathcal T$
from Theorem 1 or 2. In the examples presented below, we have used
Theorem 1 above and the analytic control of $L(s,\pi\times\pi')$
achieved by the so-called Langlands-Shahidi method of analyzing
the $L$-functions through the Fourier coefficients of Eisenstein
series \cite{S3}. Currently this requires us to take $k$ to be a
number field. The {\it functional equation}
$L(s,\pi\times\pi')=\varepsilon(s,\pi\times\pi')L(1-s,\tilde\pi\times\tilde\pi')$
has been proved in wide generality by Shahidi \cite{S1}. The {\it
boundedness in vertical strips} has been proved in close to the
same generality by Gelbart and Shahidi \cite{GS}. As for the
entire continuation of $L(s,\pi\times\pi')$, a moments thought
will tell you that one should not always expect a cuspidal
representation of $\H(\AA)$ to necessarily lift to a cuspidal
representation of $\glna$. Hence it is unreasonable to expect all
$L(s,\pi\times\pi')$ to be entire. We had previously understood
how to work around this difficulty from the point of view of
integral representations by again using a highly ramified twist.
Kim realized that  one could also control the entirety of these
twisted $L$-functions in the context of the Langlands-Shahidi
method by using a highly ramified twist. We will return to this
below. Thus in a fairly general context one has that
$L(s,\pi\times\pi')$ is {\it entire} for $\pi'$ in a suitably
modified twisting set $\mathcal T'$.


3. {\it Application of the Converse Theorem}: Once we have that
   $L(s,\pi\times\pi')$ is nice for a suitable twisting set $\mathcal
   T'$ then from the equalities
\[
L(s,\Pi\times\pi')=L(s,\pi\times\pi') \quad\quad\text{and}\quad\quad
\varepsilon(s,\Pi\times\pi')=\varepsilon(s,\pi\times\pi')
\]
we see that the $L(s,\Pi\times\pi')$ are nice and then we can apply
our Converse Theorems to conclude that $\Pi$ is either cuspidal
automorphic or at least that there is an automorphic $\Pi'$ such that
$\Pi_v=\Pi'_v$ at almost all places. This then effects the (possibly
weak) automorphic lift of $\pi$ to $\Pi$ or $\Pi'$.


4. {\it Highly ramified twists}: As we have indicated above, there are
both local and global problems that can be finessed by an appropriate
use of a highly ramified twist.
This is based on the following simple observation.


{\bf Observation.} \it Let $\Pi$ be as in Theorem 1 or 2. Suppose
that $\eta$ is a fixed  character of
$k^\times\bs\AA^\times$. Suppose that
$L(s,\Pi\times{\pi'})$ is nice for all ${\pi'}\in \mathcal T'=
\mathcal T\otimes\eta$,
where $\mathcal T$ is either of the twisting sets of Theorem 1  or
2. Then $\Pi$ is quasi-automorphic as in those theorems. \rm


The only thing to observe
 is that if ${\pi'}\in \mathcal T$ then
$L(s,\Pi\times({\pi'}\otimes\eta))=L(s,(\Pi\otimes\eta)\times{\pi'})$ so
that applying the Converse Theorem for $\Pi$ with twisting set $\mathcal
T\otimes\eta $ is equivalent to applying the Converse Theorem for
$\Pi\otimes\eta$ with the twisting set $\mathcal T$. So, by either Theorem
1 or 2, whichever is appropriate, $\Pi\otimes\eta$ is
quasi-automorphic and hence $\Pi$ is as well.

If we now begin with $\pi$ automorphic on $\H(\AA)$, we will take $T$
to be the set of finite places where $\pi_v$ is ramified. For applying
Theorem 1 we want $S=T$ and for Theorem 2 we would want $S\cap
T=\emptyset$. We will now take $\eta$ to be highly ramified at all
places $v\in T$, so that at $v\in T$ our twisting representations are all
locally of the form (unramified principal series)$\otimes$(highly
ramified character).

In order to finesse the lack of knowledge of an appropriate local lift,
we  need to know the following two local facts about the local
theory of $L$-functions for $\H$.


{\bf  Multiplicativity of \boldmath$\gamma$-factors.} \it If $\pi'_v=\Ind(\pi'_{1,v}\otimes\pi'_{2,v})$, with
$\pi'_{i,v}$ and irreducible admissible representation of $\gl_{r_i}(k_v)$, then we have
$\gamma(s,\pi_v\times\pi'_v,\psi_v)=\gamma(s,\pi_v\times\pi'_{1,v},\psi_v) \gamma(s,\pi_v\times\pi'_{2,v},\psi_v).
$ \rm


{\bf  Stability of \boldmath$\gamma$-factors.} \it
 If $\pi_{1,v}$ and $\pi_{2,v}$ are two irreducible
admissible representations of $\H(k_v)$ with the same central
 character, then for every sufficiently
highly ramified character $\eta_v$ of $GL_1(k_v)$ we have
$\gamma(s,\pi_{1,v}\times\eta_v,\psi_v)=\gamma(s,\pi_{2,v}\times\eta_v,\psi_v).
$\rm


Both of these facts are known for
$\gln$, the multiplicativity being found in \cite{JPSS} and the stability
in \cite{JS'}. Multiplicativity in a fairly wide generality useful for
applications has been established by Shahidi \cite{S4}. Stability is
in a more primitive state at the moment, but Shahidi has begun to
establish the necessary results in a general context in \cite{S2}.

To utilize these local results, what one now does is the
following. At the places where $\pi_v$ is ramified, choose $\Pi_v$
to be arbitrary, except that it should have the same central
character as $\pi_v$. This is both to guarantee that the central
character of $\Pi$ is the same as that of $\pi$ and hence
automorphic and to guarantee that the stable forms of the
$\gamma$-factors for $\pi_v$ and $\Pi_v$ agree. Now form
$\Pi=\otimes'\Pi_v$. Choose our character $\eta$ so that at the
places $v\in T$ we have that the $L$- and $\gamma$-factors for
both $\pi_v\otimes\eta_v$ and $\Pi_v\otimes\eta_v$ are in their
stable form and agree. We then twist by $\mathcal T'=\mathcal
T\otimes \eta$ for this {\it fixed} character $\eta$. If
$\pi'\in\mathcal T'$, then for $v\in T$, $\pi'_v$ is of the form
$\pi'_v=\Ind(|\ |^{s_1}\otimes\cdots\otimes |\
|^{s_m})\otimes\eta_v$. So  at the places $v\in T$, applying both
multilplcativity and stability, we have
\begin{align*}
\gamma(s,\pi_v\times\pi'_v,\psi_v)
&=\prod \gamma(s+s_i,\pi_v\otimes\eta_v,\psi_v)\\
&=\prod \gamma(s+s_i,\Pi_v\otimes\eta_v,\psi_v)
=\gamma(s,\Pi_v\times\pi'_v,\psi_v)
\end{align*}
from which one deduces a similar equality for the $L$- and
$\varepsilon$-factors. From this it will then follow that globally
we will have $L(s,\pi\times\pi')=L(s,\Pi\times\pi')$ for all
$\pi'\in\mathcal T'$ with similar equalities for the
$\varepsilon$-factors. This then completes Step 1.

To complete our use of the highly ramified twist, we must return
to the question of whether $L(s,\pi\times\pi')$ can be made
entire. In analysing $L$-functions via the Langlands-Shahidi
method, the poles of the $L$-function are controlled by those of
an Eisenstein series. In general, the inducing data for the
Eisenstein series must satisfy a type of self-contragredience for
there to be poles. The important observation of Kim is that one
can use a highly ramified twist to destroy this
self-contragredience at one place, which suffices, and hence
eliminate poles. The precise condition will depend on the
individual construction. A more detailed explanation of this can
be found in Shahidi's article \cite{S3}. This completes Step 2
above.

\section{New examples of functoriality} \label{section 4}
\setzero\vskip-5mm \hspace{5mm }

Now take $k$ to be a number field. There has been
much progress recently in utilizing the method described above to
establish global liftings from split groups $\H$ over $k$ to an
appropriate $\gln$. Among them are the following.

1. {\it Classical groups}. Take $\H$ to be a split classical group over
$k$, more  specifically, the split form of either $\SO_{2n+1}$,
   $\Sp_{2n}$, or $\SO_{2n}$. The the $L$-groups $\LH$
   are then $\Sp_{2n}(\CC)$, $\SO_{2n+1}(\CC)$, or $\SO_{2n}(\CC)$ and
   there are natural embeddings into the  general linear
   group $\gl_{2n}(\CC)$, $\gl_{2n+1}(\CC)$, or $\gl_{2n}(\CC)$
respectively.
Associated to each there should be a
lifting of admissible or automorphic representations from
$\H(\AA)$ to the appropriate $\gl_N(\AA)$. The first lifting that
resulted from the combination of the Converse Theorem and the
Langlands-Shahidi method of controlling automorphic $L$-functions was
the weak lift for generic cuspidal representations
from $\SO_{2n+1}$ to $\gl_{2n}$ over a number field $k$ obtained with
Kim and Shahidi \cite{CKPSS}.
We can now extend this to the following result.


{\bf Theorem.} \cite{CKPSS, CKPSS2}  \it Let  $\H$ be a
split classical group over $k$ as above and $\pi$ a globally generic
cuspidal representation of $\H(\AA)$. Then there exists an automorphic
representation $\Pi$ of $\gl_N(\AA)$ for the appropriate $N$ such that
$\Pi_v$ is the local Langlands lift of $\pi_v$ for all archimedean
places $v$ and almost all non-archimedean places $v$ where $\pi_v$ is
unramified. \rm


In these examples the local Langlands correspondence is not
understood at the places $v$ where $\pi_v$ is ramified and so we
must use the technique of multiplicativity and stability of the
local $\gamma$-factors as outlined in Section 3. Multiplicativity
has been established in generality by Shahidi \cite{S4} and in our
first paper \cite{CKPSS} we relied on the stability of
$\gamma$-factors for $\SO_{2n+1}$ from \cite{CPS3}. Recently
Shahidi has established an expression for his local coefficients
as Mellin transforms of Bessel functions in some generality, and
in particular in the cases at hand one can combine this with the
results of \cite{CPS3} to obtain the  necessary stability in the
other cases, leading to the extension of the lifting to the other
split classical groups \cite{CKPSS2}.


2. {\it Tensor products}. Let  $H=\gl_m\times \gl_n$.
Then $\LH=\gl_m(\CC)\times\gln(\CC)$. Then there is a natural
   simple tensor product map from $\gl_m(\CC)\times
   \gl_n(\CC)$ to  $\gl_{mn}(\CC)$. The associated functoriality from
$\gl_n\times\gl_m$ to $\gl_{mn}$ is the {\it tensor product
   lifting}. Now  the associated local lifting
   is understood in principle since the local Langlands
conjecture for $\gl_n$ has been solved.
   The question of  global functoriality has been
   recently solved in the cases of $\gl_2\times \gl_2$ to $ \gl_4$
   by Ramakrishnan \cite{R} and $\gl_2 \times \gl_3$ to  $\gl_6$ by Kim and
   Shahidi \cite{KS1, KS2}.


{\bf Theorem.} \cite{R, KS1} \it Let $\pi_1$ be a cuspidal
representation of $\gl_2(\AA)$
and $\pi_2$ a cuspidal representation of $\gl_2(\AA)$ (respectively
$\gl_3(\AA)$). Then there is an automorphic representation $\Pi$ of
$\gl_4(\AA)$ (respectively $\gl_6(\AA)$) such that $\Pi_v$ is the
local tensor product lift of $\pi_{1,v}\times\pi_{2,v}$ at all places
$v$.\rm


In both cases the authors are able to characterize when the lift is cuspidal.

In the case of Ramakrishnan \cite{R}
   $\pi=\pi_1\times\pi_2$ with each $\pi_i$ cuspidal representation
   of $\gl_2(\AA)$ and $\Pi$ is to be an automorphic representation of
   $\gl_4(\AA)$. To apply the Converse Theorem
   Ramakrishnan needs to control the analytic properties of
   $L(s,\Pi\times\pi')$ for $\pi'$ cuspidal representations of
   $\gl_1(\AA)$ and $\gl_2(\AA)$, that is, the Rankin triple product
   $L$-functions $L(s,\Pi\times\pi')=L(s,\pi_1\times\pi_2\times\pi')$.
This he was able to do using a combination of results on the integral
representation for this $L$-function due to
Garrett, Rallis and Piatetski-Shapiro, and Ikeda
and the work of Shahidi on the
Langlands-Shahidi method.

In the case of Kim and Shahidi \cite{KS1, KS2}
 $\pi_2$ is a cuspidal representation of $\gl_3(\AA)$.
Since the lifted representation $\Pi$ is to be an
automorphic representation of $\gl_6(\AA)$, to apply the Converse
Theorem  they must control the
analytic properties of $L(s,\Pi\times\pi')=L(s,\pi_1\times\pi_2\times\pi')$
where now $\pi'$ must run over appropriate cuspidal representations of
$\gl_m(\AA)$ with $m=1,2,3,4$. The control of these triple products is
an application  of the Langlands-Shahidi method of analysing
$L$-functions and involves coefficients of Eisenstein series on $\gl_5$,
$\Spin_{10}$, and simply connected $\Ex_6$ and $\Ex_7$ \cite{KS1,S3}. We should
note that even though the complete local lifting theory is understood,
they still use a highly ramified twist to control the global
properties of the $L$-functions involved. They then show that their
lifting is correct at all local places by using a base
change argument.


3. {\it Symmetric powers}. Now take  $H=\gl_2$, so $\LH=\gl_2(\CC)$.
For each $n\geq 1$ there is the
   natural symmetric $n$-th power  map $sym^n: \gl_2(\CC)\rightarrow
   \gl_{n+1}(\CC)$. The associated functoriality is the
   {\it symmetric power lifting} from representations of $\gl_2$ to
   representations of $\gl_{n+1}$. Once again
   the local symmetric powers liftings are understood in principle
   thanks to the solution of the local Langlands conjecture for
   $\gl_n$. The global symmetric square lifting,
so $\gl_2$ to $\gl_3$, is an old
   theorem of Gelbart and Jacquet. Recently, Kim and Shahidi have
   shown the existence of the global symmetric cube lifting from $\gl_2$ to
   $\gl_4$ \cite{KS1} and then Kim followed with the global symmetric
   fourth power lifting from $\gl_2$ to $\gl_5$ \cite{K}.


{\bf Theorem.} \cite{KS1,K} \it Let $\pi$ be a cuspidal
automorphic representation
of $\gl_2(\AA)$. Then there exists an automorphic representation $\Pi$
of $\gl_4(\AA)$ (resp. $\gl_5(\AA)$) such that $\Pi_v$ is the
local symmetric cube (resp. symmetric fourth power) lifting of
$\pi_v$. \rm


In either case, Kim and Shahidi have been able to give a very
interesting characterization of when the image is in fact cuspidal
\cite{KS1,KS2}.

 The original symmetric square lifting of
   Gelbart and Jacquet indeed used the converse theorem for $\gl_3$.
For Kim and Shahidi, the symmetric cube
was deduced from the functorial $\gl_2\times \gl_3$ tensor product lift above
\cite{KS1, KS2} and did not require a new use of the Converse Theorem.
For the symmetric fourth power lift, Kim first used the Converse
Theorem to establish the {\it exterior square} lift from $\gl_4$ to $\gl_6$
by the method outlined above and then combined this with the symmetric
cube lift to deduce the symmetric fourth power lift \cite{K}.

\section{Applications} \label{section 5}
\setzero\vskip-5mm \hspace{5mm }

These new examples of functoriality have already had many applications.
We will discuss the primary applications in parallel with
our presentation of the examples. $k$ remains a number field.


{\it 1. Classical groups}: The applications so far of the lifting
from classical groups to $\gln$ have been ``internal'' to the
theory of automorphic forms. In the case of the lifting from
$\SO_{2n+1}$ to $\gl_{2n}$, once the weak lift is established,
then the theory of Ginzburg, Rallis, and Soudry \cite{GRS} allows
one to show that this weak lift is indeed a strong lift in the
sense that the local components $\Pi_v$ at those $v\in S$ are
completely determined and  to completely characterize the image
locally and globally. This will be true for the liftings from the
other classical groups as well. Once one knows that these lifts
are rigid, then one can begin to define and analyse the  local
lift for ramified representations by setting the lift of $\pi_v$
to be the $\Pi_v$ determined by the global lift. This is the
content of the papers of Jiang and Soudry \cite{JS1,JS2} for the
case of $\H=SO_{2n+1}$. In essence they show that this local lift
satisfies the relations on $L$-functions that one expects from
functoriality and then deduce the {\it local Langlands conjecture
for $\SO_{2n+1}$} from that for $\gl_{2n}$. We refer to their
papers for more detail and precise statements.

{\it 2. Tensor product lifts}: Ramakrishnan's original motivation for
establishing the tensor product lifting from $\gl_2\times\gl_2$
to $\gl_4$ was to prove the multiplicity one conjecture for $\mbox{\upshape
SL}_2$ of Langlands and Labesse.


{\bf Theorem.} \cite{R} \it In the spectral decomposition
\[
L^2_{cusp}(\mbox{\upshape SL}_2(k)\backslash \mbox{\upshape
SL}_2(\mathbb A))=\bigoplus \ m_\pi \pi
\]
into irreducible cuspidal representations, the multiplicities $m_\pi$
are at most one.\rm


This was previously known to be true for $\gl_n$ and false for
$\mbox{\upshape SL}_n$ for $n\geq 3$. For further applications, for
example to  the Tate conjecture, see \cite{R}.

The primary application of the tensor product lifting from
$\gl_2\times\gl_3$ to $\gl_6$ of Kim and Shahidi was in the
establishment of the symmetric cube lifting and through this the
symmetric fourth power lifting, so
 the applications of the symmetric power
liftings outlined below are  applications of
this lifting as well.

{\it 3. Symmetric powers}: It was early observed that the existence of
the symmetric power liftings of $\gl_2$ to $\gl_{n+1}$
for all $n$ would imply the
Ramanujan-Petersson and Selberg conjectures for modular forms. Every
time a symmetric power lift is obtained we obtain better bounds
towards Ramanujan. The result which follows from the symmetric third
and fourth power lifts of Kim and Shahidi is the following.


{\bf Theorem.} \cite{KS2} \it Let $\pi$ be a cuspidal
representation of $\gl_2(\AA)$ such that the symmetric cube lift
of $\pi$ is again cuspidal. Let $\diag(\alpha_v,\beta_v)$ be the
Satake parameter for an unramified local component. Then
$|\alpha_v|, |\beta_v|<q_v^{1/9}$. If in addition the fourth
symmetric power lift is not cuspidal,  the full Ramanujan
conjecture is valid. \rm


The corresponding statement at infinite places, i.e., the analogue of the
Selberg conjecture on the eigenvalues of Mass forms, is also valid \cite{K}.
Estimates towards Ramanujan are a staple of improving any analytic
number theoretic estimates obtained through spectral methods. Both
the $1/9$ non-archimedean and $1/9$ archimedean estimate towards
Ramanujan above were applied in obtaining the precise form of the
exponent in our recent result  with Sarnak breaking the convexity
bound for twisted
Hilbert modular $L$-series in the conductor aspect, which in turn
was the key ingredient in our work on Hilbert's eleventh
problem for ternary quadratic forms. Similar in spirit
are the applications by Kim and Shahidi to the
hyperbolic circle problem and to  estimates on sums of
shifted Fourier coefficients \cite{KS1}.

In addition Kim and Shahidi were able to obtain results towards the
Sato-Tate conjecture.


{\bf Theorem.}~\cite{KS2} \it Let $\pi$ be a cuspidal representation of $\gl_2(\AA)$ with trivial central
character. Let $\diag(\alpha_v,\beta_v)$ be the Satake parameter for an unramified local component and let
$a_v=\alpha_v+\beta_v$. Assuming $\pi$ satisfies  the Ramanujan conjecture, there are sets $T^\pm$ of positive
lower density for which $a_v>2\cos(2\pi/11)-\epsilon$ for all $v\in T^+$ and $a_v<-2\cos(2\pi/11)+\epsilon$ for
all $v\in T^-$. [Note: $2\cos(2\pi/11)=1.68...$]\rm


Kim and Shahidi have other conditional applications of their
liftings such as  the conditional existence of Siegel
modular cusp forms of weight
$3$ (assuming Arthur's multiplicity formula for $\Sp_4$). We refer the
reader to \cite{KS1} for details on these applications and others.

\end{document}